\documentclass[11pt]{amsart}

\usepackage{epsfig, xypic, eucal, amssymb, amsthm, latexsym, amsmath, stmaryrd, enumerate, graphicx}

\newtheorem{theorem}{Theorem}[section]
\newtheorem{lemma}[theorem]{Lemma}

\newtheorem{corollary}[theorem]{Corollary}
\newtheorem{remark}[theorem]{Remark}

\newtheorem{conjecture}[theorem]{Conjecture}
\newtheorem{definition}[theorem]{Definition}

\def\av{\text{\rm av}}

\newtheorem*{te*}{Theorem}

\def\dsp{\def\baselinestretch{1.37}\large}

\begin{document}
\title{On the conjectures of Atiyah and Sutcliffe}

\author{Marcin Mazur}

\author{Bogdan V.~Petrenko}

\address{
Department of Mathematics \\
Binghamton University \\
P.O. Box 6000 \\
Binghamton, NY 13892-6000 } \email{ mazur@math.binghamton.edu}

\address{
Department of Mathematics \\ SUNY Brockport \\ 350 New Campus
Drive \\ Brockport, NY 14420
 }

\email{ bpetrenk@brockport.edu}

\begin{abstract}
Motivated by certain questions in physics, Atiyah defined a determinant function which to any
set of $n$ distinct points $x_1,\ldots, x_n$ in $\mathbb R^3$ assigns a complex number $D(x_1,\ldots, x_n)$.
In a joint work, he and Sutcliffe stated three intriguing conjectures about this determinant. They
provided compelling numerical evidence for the conjectures and an interesting physical interpretation
of the determinant. The first conjecture asserts that the determinant never vanishes, the second states that
its absolute value is at least one, and the third says that $|D(x_1,\ldots, x_n)|^{n-2}\geq \prod_{i=1}^n
|D(x_1,\ldots, x_{i-1},x_{i+1},\ldots, x_n)|$. Despite their simple formulation, these conjectures
appear to be notoriously difficult. Let $D_n$ denote the Atiyah determinant evaluated at the vertices of a
regular $n-$gon. We prove that $\lim_{n\to \infty} \frac{\ln D_n}{n^2}=
\frac{7\zeta(3)}{2\pi^2}-\frac{\ln 2}{2}=0.07970479...$ and establish the second conjecture in this case.
Furthermore, we prove the second conjecture for vertices of a convex quadrilateral and the third conjecture
for vertices of an inscribed quadrilateral.

\vspace{4mm}


\vspace{3mm} \noindent {\bf Keywords:} Atiyah-Sutcliffe conjecture, Atiyah determinant, configuration space.
\end{abstract}
\maketitle

\dsp
\section{Introduction}
In late 1990's, Berry and Robbins \cite{berry}, motivated by certain problems in quantum physics,
asked an interesting geometric question which can be reformulated as follows: given a positive integer $n$, is there a continuous
map which to any $n$ pairwise distinct points $x_1,\ldots ,x_n$ in $\mathbb R^3$
assigns $n$ points $p_1(x_1,\ldots, x_n),\ldots, p_n(x_1,\ldots, x_n)$
in the complex projective space $\mathbb P\mathbb C^{n-1}$ in such a way that

\begin{itemize}
\item  the points $p_1(x_1,\ldots, x_n), \ldots, p_n(x_1,\ldots, x_n)$ are not
contained in a linear subspace;

\item $p_k(x_{\sigma(1)},\ldots, x_{\sigma(n)})=p_{\sigma(k)}(x_1,\ldots, x_n)$
for any $k\in\{1,\ldots, n\}$ and any permutation $\sigma$ of $\{1,\ldots, n\}$?
\end{itemize}

The question has been answered in the positive by Atiyah in \cite{atiyah1}. In the same work Atiyah
observed that a more elegant (and more desirable) solution could be given if a certain determinant assigned to
any $n$ distinct points in $\mathbb R^3$ does not vanish. This determinant has been refined in \cite{atiyah2},
where some numerical evidence supporting its conjectural non-vanishing is given. Further refinements and
generalizations of the conjecture together with compelling numerical evidence were presented by Atiyah and
Sutcliffe in \cite{atiyah3}.  In that paper
the authors construct a determinant function with remarkable properties, which assigns to any $n$ distinct
points $x_1,\ldots, x_n$ in $\mathbb R^3$ a complex number $D(x_1,\ldots, x_n)$
(see \cite[formula (3.9)]{atiyah3}). Let us briefly outline the construction
of $D$. Denote by $(x_{i,1}, x_{i,2}, x_{i,3})$ the coordinates of $x_i$, where $i\in \{1,2,\ldots,n\}$.
For each pair $i< j$ choose two complex numbers
$z_{i,j}$ and $w_{i,j}$ such that $|z_{i,j}|^2+|w_{i,j}|^2=1$ and

\begin{equation}\label{det}
\frac{z_{i,j}}{w_{i,j}}=\frac{(x_{j,1}-x_{i,1})+(x_{j,2}-x_{i,2})\sqrt{-1}}
{\sqrt{\sum_{k=1}^3(x_{j,k}-x_{i,k})^2}-(x_{j,3}-x_{i,3})},
\end{equation}
with the convention that $w_{i,j}=0$ when the denominator of the right hand side of (\ref{det}) vanishes.
When $i>j$, define $z_{i,j}=-\overline{w}_{j,i}$ and $w_{i,j}=\overline{z}_{j,i}$. Define $a_{i,j}$ as
the coefficient at $t_1^{j-1} t_2^{n-j}$
of the polynomial $f_i(t_1,t_2)=\prod_{k\neq i}(z_{i,k}t_1-w_{i,k}t_2)$. The  Atiyah determinant $D(x_1,\ldots, x_n)$ is
defined as the determinant
of the matrix $(a_{i,j})$. In \cite{atiyah3} it has been proved that $D(x_1,\ldots, x_n)$ is independent of all
the choices made in the course of its definition.
Moreover, this determinant is invariant under the orientation preserving similitudes of $\mathbb R^3$ and becomes
its own conjugate under the orientation
reversing similitudes.
In \cite{atiyah3} the authors stated the following three conjectures about $D(x_1,\ldots, x_n)$.

\begin{conjecture}\label{conj1}
$D(x_1,\ldots, x_n)\neq 0$ for all $x_1,\ldots, x_n$.
\end{conjecture}

\begin{conjecture}\label{conj2}
$|D(x_1,\ldots, x_n)|\geq 1$ for all $x_1,\ldots, x_n$.
\end{conjecture}

\begin{conjecture}\label{conj3} For all $x_1,\ldots, x_n$ we have
\[|D(x_1,\ldots, x_n)|^{n-2}\geq \prod_{i=1}^n |D(x_1,\ldots, x_{i-1},x_{i+1},\ldots, x_n)|.\]
\end{conjecture}

It is easy to see that the conjectures are stated in order of increasing strength. All three conjectures have
been verified by Atiyah for $n=3$. In \cite{atiyah3} a compelling numerical evidence is given in support
of all three conjectures. In addition, the authors provide a very interesting physical interpretation
of Atiyah determinant and discuss further generalization
of the conjectures. Conjecture~\ref{conj1} has been proved for $n=4$ by Eastwood and Norbury
\cite{east}. In addition, Conjecture~\ref{conj1} has been proved for some configurations of points of
arbitrarily large size in \cite{dok}.
We are not aware of any other results concerning these conjectures.

In the first part of our paper we obtain an explicit formula for the value $D_n$ of Atiyah determinant
at vertices of any regular $n$-gon (see Theorem~\ref{formul}).  Using this formula we prove in
Theorem~\ref{asymp} that $\lim_{n\to \infty} \frac{\ln D_n}{n^2}=
\frac{7\zeta(3)}{2\pi^2}-\frac{\ln 2}{2}\approx 0.07970479$ and
confirm Conjecture~\ref{conj2} in this case.
Note that Conjecture~\ref{conj1} in this case follows from the results of \cite{dok}.
In the second part of the paper, building on the
work of Eastwood and Norbury \cite{east}, we investigate Atiyah determinant when $n=4$. In Theorem~\ref{convex} we
prove Conjecture~\ref{conj2} for vertices of any convex quadrilateral, and in Theorem~\ref{inscribed} we confirm
Conjecture~\ref{conj3} for inscribed quadrilaterals. In the course of proving these results, we are led
to some intriguing results and conjectures about tetrahedra and quadrilaterals for which we have compelling
numerical evidence (see Conjectures~\ref{conj4},~\ref{conj5}, and~\ref{conj6}).

\section{Regular n-gon}
Suppose that the points $x_1, \ldots, x_n$ are on a circle. Recall that Atiyah determinant
is invariant under orientation preserving similitudes of $\mathbb R^3$. Therefore, in order to compute
$D(x_1,\ldots, x_n)$, we may assume that $x_{i,3}=0$ for all $i$ and
$x_{i,1}+x_{i,2}\sqrt{-1}=e^{u_k\sqrt{-1}}$, where $0<u_1<\ldots<u_n\leq 2\pi$.
Set $a_r=e^{u_r\sqrt{-1}/2}$ for
$r=1,\ldots, n$. A straightforward calculation
confirms that we can take $z_{i,j}=a_i/\sqrt{2}$, $w_{i,j}=\sqrt{-1}a_j^{-1}/\sqrt{2}$ for $1\leq i<j\leq n$
in the computation of $D$. Define $g_r(z)=\prod_{s<r}(z+a_ra_s)\prod_{s>r}(z-a_ra_s)$.
Then
\begin{equation}\label{poly4}
f_r(t_1,t_2)= \sqrt{2}^{1-n}(\sqrt{-1})^{r-1}(-1)^{n-r}\left(\prod_{k\neq r}a_k\right)^{-1}t_1^{n-1}g_r\left(
\sqrt{-1}\frac{t_2}{t_1}\right).
\end{equation}

It follows easily from the last formula that Conjecture~\ref{conj1} for the points $x_1,\ldots, x_n$ is
equivalent to $\mathbb C$-linear independence of the polynomials $g_r(z)$, $r=1,\ldots, n$. Indeed, by
(\ref{poly4}), the $\mathbb C$-linear independence of these polynomials is equivalent
to the $\mathbb C$-linear independence of the polynomials $f_r(t_1, t_2)$, $r=1,\ldots, n$.
In turn, the $\mathbb C$-linear independence of the latter sequence of polynomials is, by definition,
equivalent to the non-vanishing of the determinant $D(x_1,\ldots, x_n)$.

We specialize now to the case when the points $x_1,\ldots, x_n$ are vertices of a regular $n$-gon.
In other words, we assume that $u_k= 2\pi k/n$, $k=1,\ldots, n$.
Define $g(z)=\prod_{k=1}^{n-1}(z-w^k)$, where $w=e^{\pi \sqrt{-1}/n}$. Then $a_r=w^r$ and
$g_r(z)=w^{2r(n-1)}g(w^{-2r}z)$, $r=1,\ldots, n$.

\begin{lemma}\label{nonzero}
Suppose that $a\neq 0$ and $a^k\neq 1$ for all $k$ such that  $1\leq k < n$.
Let $h(z)=\prod_{k=1}^n(z-a^k)=z^n-\sum_{k=0}^{n-1}b_k z^k$. Then
\begin{equation}\label{vander}
 b_k=a^{n-k}\prod_{l\neq k, 0\leq l<n}(a^n-a^l)\prod_{l\neq k, 0\leq l< n}(a^k-a^l)^{-1}.
 \end{equation}
\end{lemma}
\begin{proof}
Let $V(y_1,\ldots, y_n)$ be the Vandermonde matrix, i.e. the  $n\times n$ matrix whose $(k,l)$-entry
is $y_k^{l-1}$. Recall that the determinant of this matrix is given by
\begin{equation}\label{vander1}
\det V(y_1,\ldots, y_n)=\prod_{1\leq i<j\leq n}(y_j-y_i).
\end{equation}
The equalities $h(a^m)=0$, $m=1,\ldots, n$ translate
into a system of $n$ linear equations for the coefficients $b_k$:
\[ a^{mn}=\sum_{k=0}^{n-1}a^{mk}b_k,\ m=1,\ldots, n.\]
Using Cramer's rule and formula (\ref{vander1}), it is a straightforward computation to get (\ref{vander}).
\end{proof}

\begin{corollary}
Conjecture~\ref{conj1} is true for vertices of a regular $n$-gon.
\end{corollary}

\begin{proof}
It suffices to prove that the polynomials $g(w^{-2r}z)$, $r=1,\ldots,n$, are linearly independent
over $\mathbb C$. Write $g(z)=z^{n-1}+\sum_{t=0}^{n-2}b_tz^t$. By Lemma~\ref{nonzero},
all the coefficients $b_t$ are non-zero. Thus the equality  $\sum_{r=1}^n x_rg(w^{-2r}z)=0$
is equivalent to the system of $n$ linear equations:
\[ \sum_{r=1}^{n} x_rw^{-2rt}=0, \ t=0,1,\ldots, n-1.\]
Therefore the polynomial $\sum_{r=1}^{n}x_rz^{r-1}$ of degree $n-1$ has $n$ distinct
roots $w^{-2t}$, $t=0,1,\ldots, n-1$. It follows that this polynomial is $0$, i.e. $x_1=\ldots = x_n=0$.
This establishes the linear independence of the polynomials $g(w^{-2r}z)$, $r=1,\ldots, n$.
\end{proof}

We are now ready to compute $D(x_1,\ldots, x_n)$.

\begin{theorem}\label{formul}
Let $x_1,\ldots, x_n$ be vertices of a regular $n$-gon. Then
\begin{equation}\label{formula5}
|D(x_1, \ldots, x_n)|=n^{n/2}2^{n(1-n)/2}\prod_{1\leq k\leq n/2}\left(\cot\frac{\pi k}{2n}\right)^{n-2k}.
\end{equation}
\end{theorem}

\begin{proof}
In order to carry out the computation of $D(x_1,\ldots, x_n)$ note that (\ref{poly4})
yields
\[ f_r(t_1,t_2)=-\sqrt{2}^{1-n}(\sqrt{-1})^{n-r}w^{-r}t_1^{n-1}g\left(\sqrt{-1}w^{-2r}\frac{t_2}{t_1}\right)
\]
(we used the equality $\left(\prod_{k\neq r}a_k\right)^{-1}= \left(\prod_{k\neq r}w^k\right)^{-1}=w^r(-\sqrt{-1})^{n+1}$).
If $g(z)=z^{n-1}-\sum_{t=0}^{n-2}b_tz^t$ then the entries $a_{r,j}$ of the matrix defining
$D(x_1, \ldots, x_n)$ are given by
\[a_{r,j}=\sqrt{2}^{1-n}(\sqrt{-1})^{n-r}w^{-r}(\sqrt{-1}w^{-2r})^{n-j}b_{n-j},
\]
where we set $b_{n-1}=-1$. Thus
\[|D(x_1, \ldots, x_n)|=\sqrt{2}^{n(1-n)}\left|\det V(
w^{-2\cdot 1},w^{-2\cdot 2},\ldots, w^{-2\cdot n})\prod_{i=0}^{n-2}b_i\right|.
\]
A straightforward computation, using (\ref{vander1}) and the identity $\prod_{s=1}^{2n-1}(1-w^s)=2n$, yields
\[\left|\det V(w^{-2\cdot 1},w^{-2\cdot 2},\ldots, w^{-2\cdot n})\right|= n\prod_{0\leq s<t\leq n-2}
|w^{2t}-w^{2s}|.
\]
Using Lemma~\ref{nonzero}, we get that
\[\prod_{i=0}^{n-2}|b_i|=\prod_{s=1}^{n-1}|1-w^s|^{n-2}
\prod_{0\leq s<t\leq n-2}|w^{t}-w^{s}|^{-2}.
\]
Since $\overline{1-w^s}=1-w^{2n-s}$, we have the following equality:
\[ 2n=\prod_{s=1}^{2n-1}(1-w^s)=(1-w^n)\prod_{s=1}^{n-1}|1-w^s|^2=2\prod_{s=1}^{n-1}|1-w^s|^2.\]
Putting all these computations together, we arrive at the following formula:
\[|D(x_1,\ldots, x_n)|=n^{n/2}2^{n(1-n)/2}\prod_{0\leq s<t\leq n-2}\left|\frac{w^{t}+w^{s}}{w^{t}-w^{s}}\right|.
\]
A straightforward calculation, using the identity $\displaystyle \left|\frac{w^{t}+w^{s}}{w^{t}-w^{s}}\right|
=\frac{1+e^{\alpha\sqrt{-1}}}
{1-e^{\alpha\sqrt{-1}}}=\cot(\alpha/2)$ (for an appropriate $\alpha$),
yields  (\ref{formula5}).
\end{proof}

Our next goal is to confirm Conjecture~\ref{conj2} for the vertices of a regular $n$-gon. We start with some lemmas.

\begin{lemma}\label{decreasing}
The function $\displaystyle f(x)=\left(\frac{\pi}{4}-x\right)\ln\cot x$ is decreasing on $(0,\pi/4)$.
\end{lemma}
\begin{proof}
We have $\displaystyle f'(x)=-\ln\cot x-\left(\frac{\pi}{2}-2x\right)\csc 2x$. It suffices to show that
$f'(x)<0$ on $(0,\pi/4)$. This is equivalent to showing that $g(x):=2x-\sin 2x \ln\cot x<\pi/2$.
Now $g'(x)= 4-2\cos 2x\ln\cot x$ and $g''(x)=4\sin 2x\ln\cot x+4\cot 2x$. It is clear that $g''(x)>0$
on $(0,\pi/4)$. Thus $g$ is concave up on $(0,\pi/4)$ so its largest value on the interval $[0,\pi/4]$
is attained at one of the ends. Note that $\lim_{x\to 0^+}g(x)=0$ and $g(\pi/4)=\pi/2$, which proves our
claim.
\end{proof}

\begin{lemma}\label{integral} Let $\zeta(x)$ be the Riemann's zeta function. Then
\[ \int_{0}^{\pi/4}\left(\frac{\pi}{4}-x\right)\ln\cot xdx=\frac{7}{16}\zeta(3)=0.5258998951... \ \  .
\]
\end{lemma}
\begin{proof}
Integration by parts followed by a simple substitution yield
\[
\int_{0}^{\pi/4}\left(\frac{\pi}{4}-x\right)\ln\cot xdx = \frac{\pi}{8}\int_{0}^{\pi/2}
x\csc x dx-\frac{1}{8}\int_{0}^{\pi/2}
x^2\csc x dx.
\]
It turns out that both integrals on the right can be found in the literature. We have found them first
in the wonderful monograph \cite{constants}, where on pages $56-57$ the following formulas are
given (without proof):
\begin{equation}\label{sec}
\int_{0}^{\pi/2}x\csc x dx= 2G
\end{equation}
and
\begin{equation}\label{xsec}
\int_{0}^{\pi/2}x^2\csc x dx=2\pi G-\frac{7}{2}\zeta(3),
\end{equation}
where $G$ is Catalan's constant. Both formulas are proved in \cite{bradley} and (\ref{xsec}) is proved in
\cite{integr}. It is clear now that the lemma follows from (\ref{sec}) and (\ref{xsec}).
\end{proof}

\begin{lemma}\label{summ}
Let $\displaystyle B=\frac{7\zeta(3)}{2\pi^2}=0.42627839... $.
Then
\begin{equation}\label{sum}
e^{Bn-(1-\frac{1}{n})\ln n-(1-\ln(\pi/2))} \leq \prod_{1\leq k\leq n/2}\left(\cot\frac{\pi k}{2n}\right)^{1-\frac{2k}{n}}\leq e^{Bn}.
\end{equation}
\end{lemma}
\begin{proof}

Note that
\begin{equation}\label{a}
\sum_{1\leq k\leq n/2}\left(1-\frac{2k}{n}\right)\ln\cot\frac{\pi k}{2n}=
\frac{8n}{\pi^2}\sum_{1\leq k\leq n/2}f\left(\frac{\pi k}{2n}\right)\left(\frac{\pi (k+1)}{2n}-
\frac{\pi k}{2n}\right),
\end{equation}
where $\displaystyle f(x)=\left(\frac{\pi}{4}-x\right)\ln\cot x$.
The sum on the right hand side of (\ref{a}) is a Riemann sum for $f$. By Lemma~\ref{decreasing},
the function $f$ is decreasing and
non-negative on $(0,\pi/4)$. Thus
\[ \int_{\pi/2n}^{\pi/4}f(x)dx
\leq \sum_{1\leq k\leq n/2}f\left(\frac{\pi k}{2n}\right)\left(\frac{\pi (k+1)}{2n}-
\frac{\pi k}{2n}\right)\leq \int_{0}^{\pi/4}f(x)dx.
\]
It follows from Lemma~\ref{integral} that
\[Bn-\frac{8n}{\pi^2}\int_{0}^{\pi/2n}f(x)dx\leq
\sum_{1\leq k\leq n/2}\left(1-\frac{2k}{n}\right)\ln\cot\frac{\pi k}{2n}\leq Bn
\]
Using the inequality $x<\tan x$ , we see that
\[\int_{0}^{\epsilon}f(x)dx\leq \int_{0}^{\epsilon}\left(x-\frac{\pi}{4}\right)\ln x=\frac{\epsilon}{4}
(2\epsilon-\pi)\ln \epsilon+\frac{\epsilon}{4}(\pi-\epsilon). \]
For $\epsilon=\pi/2n$, we get
\[\frac{8n}{\pi^2}\int_{0}^{\pi/2n}f(x)dx\leq (1-\frac{1}{n})\ln n+1-\ln(\pi/2). \]
Thus
\[Bn-\left(1-\frac{1}{n}\right)\ln n-1+\ln(\pi/2) \leq \sum_{1\leq k\leq n/2}\left(1-\frac{2k}{n}\right)
\ln\cot\frac{\pi k}{2n}\leq Bn.
\]
Exponentiation of all sides yields (\ref{sum}).
\end{proof}

Let us note the following interesting corollary.

\begin{theorem}\label{lim}
\[ \lim_{n\to \infty} \prod_{1\leq k\leq n/2}\left(\cot\frac{\pi k}{2n}\right)^{\frac{n-2k}{n^2}}=
e^{\frac{7\zeta(3)}{2\pi^2}}.\]
\end{theorem}

We can now state and prove the main result of this section.

\begin{theorem}\label{asymp}
Let $D_n=|D(x_1,\ldots, x_n)|$, where $x_1,\ldots, x_n$ are the vertices of a regular $n$-gon. Then
\begin{equation}\label{limit}
\lim_{n\to \infty} \frac{\ln D_n}{n^2}=\frac{7\zeta(3)}{2\pi^2}-\frac{\ln 2}{2}=0.07970479...
\end{equation}
and $D_n>1$ for all $n\geq 3$.

\end{theorem}
\begin{proof}
Formula (\ref{limit}) is a straightforward consequence of Theorems~\ref{formul} and~\ref{lim}.
It remains to prove that $D_n > 1$.
By (\ref{sum}) and Theorem~\ref{formul}, we have
\[ \frac{\ln D_n}{n}\geq \left(\frac{7\zeta(3)}{2\pi^2}-\frac{\ln 2}{2}\right)n -\left(\frac{1}{2}-
\frac{1}{n}\right)\ln n+\ln\left(\frac{\pi}{\sqrt{2}}\right)-1\geq \]
\[\geq 0.0797\cdot n -\left(\frac{1}{2}-
\frac{1}{n}\right)\ln n -0.2019.
\]
It is a simple calculus exercise to see that the rightmost expression is positive and increasing with $n$
for $n\geq 20$. This implies that $D_n>1$ for $n\geq 20$. For $n<20$ the inequality $D_n>1$ is verified
by a direct computation.
\end{proof}

\begin{remark}{\rm
The fact that $D_n$ grows so rapidly should not come as a surprise. Note that the numerical investigation
in \cite{atiyah3} found $D_n$ to be the maximum of $|D(x_1,\ldots, x_n)|$ among all coplanar points
$x_1,\ldots, x_n$ for $n\leq 15$. For $n\geq 16$ this is no longer true and the investigation of
\cite{atiyah3} suggests a rather intriguing pattern for the coplanar configuration with maximal $|D(x_1,\ldots, x_n)|$.
On the other hand, when $x_1, \ldots, x_n$ are collinear, we have $|D(x_1,\ldots, x_n)|=1$, so Conjecture~\ref{conj2}
is the best possible.
}

\end{remark}

\section{Four coplanar points}
Conjecture~\ref{conj1} has been confirmed for $n=4$ in \cite{east}. The main idea of that paper is to
express Atiyah determinant $D(x_1,x_2,x_3,x_4)$ as a function of the Euclidean distances $r_{i,j}=|x_i-x_j|$.
It turns out that $[|D(x_1,x_2,x_3,x_4)|\prod_{1\leq i<j\leq 4}(2r_{i,j})]^2$ is a homogeneous
polynomial of degree 12 with 4500 terms (the authors used Maple to compute the polynomial). However,
the real part of $D(x_1,x_2,x_3,x_4)\prod_{1\leq i<j\leq 4}(2r_{i,j})$
is a much more accessible polynomial, homogeneous of degree 6 with 226 terms. In order to write this polynomial
in a compact form we recall the following notation from \cite{east}. If $f$ is a polynomial in the
variables $r_{i,j}$ (where $r_{i,j}=r_{j,i}$) and $\sigma$ is a permutation of the set $\{1,2,3,4\}$ then
$f^{\sigma}$ is
obtained from $f$ by replacing $r_{i,j}$ with $r_{\sigma(i),\sigma(j)}$ for each pair $i<j$. For example,
if $f=r_{1,3}+r_{1,4}$ and $\sigma $ is the $4$-cycle $(1,2,3,4)$ then $f^{\sigma}=r_{2,4}+r_{1,2}$.
We define $\av(f)=(\sum f^{\sigma})/24$, where the sum is over all permutations of the set $\{1,2,3,4\}$.
Finally, let
\[d_3(a,b,c)=(a+b-c)(a+c-b)(b+c-a)\]
and let $V$ be the volume of the tetrahedron with
vertices $x_1, x_2, x_3, x_4$. With this notation the real part of
$D(x_1,x_2,x_3,x_4)\prod_{1\leq i<j\leq 4}(2r_{i,j})$ is given by the following formula:
\begin{eqnarray}\label{poly}
64r_{1,2}r_{1,3}r_{1,4}r_{2,3}r_{2,4}r_{3,4}-4d_3(r_{1,2}r_{3,4}, r_{1,3}r_{2,4}, r_{1,4}r_{2,3})+ \\ \nonumber
+12\av\left(r_{1,4}((r_{2,4}+r_{3,4})^2-r_{2,3}^2)d_3(r_{1,2},r_{1,3},r_{2,3})\right)+288V^2.
\end{eqnarray}
Consider the (possibly degenerate) tetrahedron with vertices  $x_1,x_2,x_3,x_4$.
{\bf In what follows, the set of indices $\{i,j,k,l\}$ will always coincide with $\{1,2,3,4\}$.} Let $\alpha_{i,j}$
be the angle $\angle x_{k}x_ix_l$ (i.e. the angle at vertex $x_i$ of the face subtended the vertex $x_j$ ).

\begin{lemma}\label{triangle}
Let $ABC$ be a triangle with sides $a=BC$, $b=AC$, and $c=AB$. Then
\[ d_3(ABC) := d_3(a,b,c)=2abc(\cos A+\cos B+\cos C-1). \]
\end{lemma}
\begin{proof}
By the law of cosines, we have $2abc\cos A=ab^2+ac^2-a^3$ and similar identities hold for the other
two angles. The conclusion of the lemma follows now easily by adding these identities.
\end{proof}

By the law of cosines,
\[(r_{i,j}+r_{j,k})^2-r_{i,k}^2=2r_{i,j}r_{j,k}(1+\cos\alpha_{j,l}).\]
Together with Lemma~\ref{triangle} this yields
\begin{eqnarray}\label{poly1}
12\av\left(r_{1,4}((r_{2,4}+r_{3,4})^2-r_{2,3}^2)d_3(r_{1,2},r_{1,3},r_{2,3})\right)= \\ \nonumber
4\left(\prod_{1\leq i<j\leq 4}r_{i,j}\right)\sum_{l=1}^{4} (3+\cos\alpha_{l,i}+\cos\alpha_{l,j}+
\cos\alpha_{l,k})(\cos\alpha_{i,l}+\cos\alpha_{j,l}+
\cos\alpha_{k,l}-1).
\end{eqnarray}

In order to get some insight into $d_3(r_{1,2}r_{3,4}, r_{1,3}r_{2,4}, r_{1,4}r_{2,3})$ we need
the following old result about tetrahedra.

\begin{lemma}\label{crelle}
Let $ABCD$ be a tetrahedron. There exists a triangle with side lengths $AB\cdot CD$,
$AC\cdot BD$, $AD\cdot BC$. For any vertex of the tetrahedron, the angles of this triangle are equal to the
angles between the circles circumscribed on the three faces of the tetrahedron sharing the chosen vertex.
\end{lemma}

For the convenience of the reader we provide a sketch of a proof.

\begin{proof}
Pick a vertex, say $D$, of the tetrahedron and let $A'$, $B'$, $C'$ be the images of $A,B,C$
respectively under the inversion $I$ in a sphere with center $D$ and radius $r=\sqrt{DA\cdot DB\cdot DC}$.
Using the fact that $I(X)I(Y)\cdot DX\cdot DY= r^2XY$, we get
\[ A'B'=AB\cdot CD, \ B'C'=AD\cdot BC, \ C'A'=AC\cdot BD.
\]
Thus $A'B'C'$ is the required triangle. Note that $I$ takes the lines $A'B'$, $A'C'$, $B'C'$ to circles
circumscribed on the three faces of $ABC$ sharing the vertex $D$. Since $I$ is conformal, the claim about
angles of $A'B'C'$ follows.
\end{proof}

\begin{remark}\label{remcrelle}
{\rm
\

\begin{enumerate}

\item The description of the angles in Lemma~\ref{crelle} may seem ambiguous, since two intersecting
circles or lines do not define a unique angle but a pair of supplementary angles.
However, given three pairs $\{\alpha_i, \pi-\alpha_i\}$, $i=1,2,3$, of supplementary
angles, none of which is $0$, there is at most one choice of $\beta_i\in \{\alpha_i, \pi-\alpha_i\}$ such that $\beta_1+
\beta_2+\beta_3=\pi$.

\item It is clear that Lemma~\ref{crelle} remains true for degenerate tetrahedra (when $A,B,C,D$ are
coplanar) except that in this case the triangle may be degenerate, which happens if and only if the
points $A,B,C,D$ are on one line or circle.

\item It is a result of Crelle \cite{crel} that the area $S$ of the triangle in Lemma~\ref{crelle}, the volume
$V$ of the tetrahedron
and the radius $R$ of the sphere circumscribed on the tetrahedron are related by the formula $S=6VR$.
\end{enumerate}
}
\end{remark}

\begin{definition}
Given any four distinct points in  $\mathbb R^3$, the triangle discussed in Lemma~\ref{crelle} and Remark~\ref{remcrelle} (2)
will be called the {\bf Crelle triangle} associated to the four points.
\end{definition}

It follows that
\begin{equation}\label{assoc}
d_3(r_{1,2}r_{3,4}, r_{1,3}r_{2,4}, r_{1,4}r_{2,3})=2
(\cos A+\cos B+\cos C-1)\prod_{1\leq i<j\leq 4}r_{i,j},
\end{equation}
where $A,B,C$ are the angles of the associated Crelle triangle. We get the following corollary.

\begin{corollary}
Let $x_1$, $x_2$, $x_3$, $x_4$ be distinct points in $\mathbb R^3$ and let $A, B, C$ be the
angles of the associated Crelle triangle. If
\begin{eqnarray}\label{formula}
\sum_{l=1}^{4} (3+\cos\alpha_{l,i}+\cos\alpha_{l,j}+
\cos\alpha_{l,k})(\cos\alpha_{i,l}+\cos\alpha_{j,l}+
\cos\alpha_{k,l}-1) \geq \\ \nonumber
\geq 2(\cos A+\cos B+\cos C-1)
\end{eqnarray}
then Conjecture~\ref{conj2} holds for $x_1$, $x_2$, $x_3$, $x_4$.
When, in addition, the points $x_1$, $x_2$, $x_3$, $x_4$ are coplanar then (\ref{formula}) is in fact
equivalent to Conjecture~\ref{conj2}.
\end{corollary}
\begin{proof}
Since (\ref{poly}) is the real part of $\displaystyle 64
D(x_1,x_2,x_3,x_4)\prod_{1\leq i<j\leq 4}r_{i,j}$, the inequality $|D(x_1,x_2,x_3,x_4)|\geq 1$ will hold if
\begin{equation}\label{formula1}
12\av\left(r_{1,4}((r_{2,4}+r_{3,4})^2-r_{2,3}^2)d_3(r_{1,2},r_{1,3},r_{2,3})\right)\geq
4d_3(r_{1,2}r_{3,4}, r_{1,3}r_{2,4}, r_{1,4}r_{2,3}).
\end{equation}
In addition, if $x_1$, $x_2$, $x_3$, $x_4$ are coplanar then $V=0$ and $D(x_1,x_2,x_3,x_4)$ is real
so $|D(x_1,x_2,x_3,x_4)|\geq 1$ is equivalent to (\ref{formula1}). To complete the proof note that
our considerations above show that (\ref{formula1}) and  (\ref{formula}) are equivalent.

\end{proof}

We do not know any explicit formulas expressing the angles of the associated Crelle triangle in terms of the angles
$\alpha_{i,j}$ in general. However, when the points $x_1,x_2,x_3,x_4$ are coplanar, such formulas are
easy to obtain using Lemma~\ref{crelle} (or rather Remark~\ref{remcrelle} (2)).

\begin{lemma}\label{angles}
Let $x_1$, $x_2$, $x_3$, $x_4$ be distinct coplanar points.

\begin{enumerate}[\rm (i)]
\item If $x_1x_2x_3x_4$ is a convex quadrilateral and $\alpha_{1,3}+\alpha_{3,1}\leq \pi$ then the associated
Crelle triangle has angles $\alpha_{2,3}-\alpha_{3,2}$, $\alpha_{2,1}-\alpha_{1,2}$, $\alpha_{1,3}+\alpha_{3,1}$.

\item If $x_4$ belongs to the triangle $x_1x_2x_3$ then the associated Crelle
triangle has angles $\alpha_{1,2}+\alpha_{2,1}$, $\alpha_{1,3}+\alpha_{3,1}$, $\alpha_{2,3}+\alpha_{3,2}$.

\end{enumerate}
\end{lemma}
\begin{proof}
The lemma follows easily from the following fact from elementary plane geometry.
Let $c_1,c_2$ be circles intersecting in 2 points $A, B$. Let $C_1\in c_1$, $C_2\in c_2$
be points on the same side of the line $AB$. The angle between $c_1$ and $c_2$ is equal to the angle
between the lines tangent to $c_1$ and $c_2$ at the point $A$. Using the result about the angle between a
tangent and a secant (Proposition 32 in Book III of the {\em Elements}) we get that the angle between
$c_1$ and $c_2$ is $|\angle AC_1B-\angle AC_2B|$.

We leave further details to the reader. Working with directed angles may simplify the argument and
Remark~\ref{remcrelle} (1) may be useful.
\end{proof}

\begin{lemma}\label{technical}
Let the function $f(u,w,x,y,z)$ be defined as follows:
\[ f(u,w,x,y,z)=\cos u+\cos w+\cos x+\cos y +\cos z-\cos(u+y+z)-\cos(x+y+z)+\]
\[+\cos(-w+y+z)+\cos(u+w)
+\cos(x+y)-\cos(u+y)-\cos(w+x)+\cos(u+x+y+z)-\]
\[\cos(-w+z)-\cos(u+w+x+y).\]
Then $f\geq 3$ for any non-negative $u,w,x,y,z$ such that
\[w\leq z,\  x+w\leq \pi,\  u+w+x+y+z\leq 2\pi,\
u+x+y+z\leq \pi,\  \text{and}\  u+y+z\leq \pi.\]
\end{lemma}
\begin{proof}
We consider first the case when $u=0$.
\[ f(0,w,x,y,z)=1+2\cos w +\cos x-\cos(y+z)+\cos(-w+y+z)-\cos(w+x)+\]
\[+\cos z+\cos(x+y)-\cos(-w+z)-\cos(w+x+y)=1+2\cos w +\cos x-\cos(y+z)+\]
\[\cos(-w+y+z)-\cos(w+x)+4\sin(w/2)\cos[(x+y+z)/2]\sin[w+x+y-z)/2].
\]
It follows that
\[ f(0,w,x,y,z)-f(0,w,x,0,y+z)=\]
\[8\sin(w/2)\cos[(x+y+z)/2]\sin(y/2)\cos[(w+x-z)/2]\geq 0.
\]
Now
\[f(0,w,x,0,y+z)=1+2\cos w+2\cos x-2\cos(w+x)=\]
\[ 3+2(\cos w+\cos x+\cos(\pi-w-x)-1)
\]
Since $w,x, \pi-w-x$ are angles of a triangle, Lemma~\ref{triangle} allows us to conclude that
\begin{equation}\label{zero}
f(0,w,x,y,z)\geq f(0,w,x,0,y+z)\geq 3.
\end{equation}

In order to handle the general case, note that
\[h(u,w,x,y,z):=\cos z+ \cos(u+w)-\cos(x+y+z)-\cos(u+w+x+y)=\]
\[4\sin[(u+w+x+y+z)/2]\sin[(x+y)/2]\cos[(u-z+w)/2],\]
and
\[ g(u,w,x,y,z):=\cos u-\cos(-w+x)+\cos(-w+y+z)-\cos(u+y)=\]
\[4\sin(y/2)\cos[(u-w+y+z)/2] \sin[(u-z+w)/2].\]

\noindent
Let $A=4\sin[(u+w+x+y+z)/2]\sin[(x+y)/2]$ and $B=4\sin(y/2)\cos[(u-w+y+z)/2] $ and $R=\sqrt{A^2+B^2}$.
Then $A\geq 0$ and $B\geq 0$ so there is $\alpha\in [0,\pi/2]$ such that
$\sin\alpha=A/R$ and $\cos\alpha=B/R$. Then
\[h(u,w,x,y,z)+g(u,w,x,y,z)=4R\sin[\alpha+(u-z+w)/2]
\]
Note now that $f(u,w,x,y,z)$ differs from $h(u,w,x,y,z)+g(u,w,x,y,z)$ only by terms which are
functions of $w,x,y$, and $u+z$. It follows that
\[f(u,w,x,y,z)-f(0,w,x,y,u+z)=8R\sin(u/2)\cos(\alpha+(w-z)/2)\geq 0,\]
since by our assumptions we have $-\pi/2\leq -z/2\leq \alpha+(w-z)/2\leq \alpha\leq \pi/2$.
Together with (\ref{zero}), this completes the proof of the lemma.

\end{proof}

In order to state our next result more efficiently we introduce the following definition.

\begin{definition}
Let $ABC$ be a triangle with sides $a=BC$, $b=AC$, and $c=AB$. Then
\[ \delta(ABC) := \frac{d_3(a,b,c)}{2abc}=\cos A+\cos B+\cos C-1. \]
\end{definition}

\begin{remark}\label{radii}
{\rm Using Heron's formula $16S^2=d_3(a,b,c)(a+b+c)$ for the area $S$ of the triangle $ABC$ and the
formulas $4S=(abc)/R=2(a+b+c)r$, where $R$ and $r$ are radii of the circumscribed and inscribed circles
respectively, we get a nice geometric interpretation of $\delta$: $\delta(ABC)=r/R$.
}
\end{remark}

\begin{remark}\label{3points}
{\rm It is not hard to see that $\displaystyle D(A,B,C)=1+\frac{\delta(ABC)}{4}$.

}
\end{remark}

We can now state the first main result of this section.

\begin{theorem}
Let $x_1x_2x_3x_4$ be a convex quadrilateral and let $ABC$ be the associated Crelle triangle. Then
\begin{equation}\label{inequality}
\delta(x_1x_2x_3)+\delta(x_1x_3x_4)+\delta(x_1x_2x_4)+\delta(x_2x_3x_4)\geq \delta(ABC).
\end{equation}
\end{theorem}
\begin{proof}
We may assume that $\alpha_{13}+\alpha_{31}\leq \pi$ (since the sum of the angles of any quadrilateral
is $2\pi$). By Lemma~\ref{angles}, we have
\[ \delta(ABC)=\cos(\alpha_{2,3}-\alpha_{3,2})+\cos(\alpha_{2,1}-\alpha_{1,2})+\cos(\alpha_{1,3}+\alpha_{3,1})-1.
\]
\noindent
Setting $u=\alpha_{3,4}$, $w=\alpha_{3,2}$, $x=\alpha_{1,2}$, $y=\alpha_{1,4}$, $z=\alpha_{2,3}$ (see the picture
below) it is
straightforward to see that
\[\delta(x_1x_2x_3)+\delta(x_1x_3x_4)+\delta(x_1x_2x_4)+\delta(x_2x_3x_4)-\delta(ABC)=f(u,w,x,y,z)-3,\]

\vspace{3mm}

\center \includegraphics{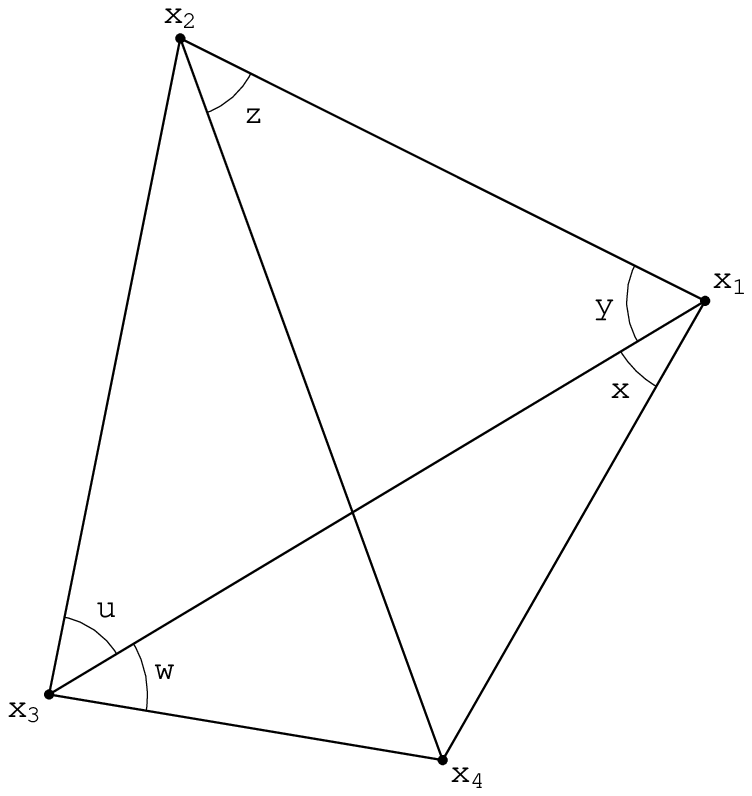}

where $f$ is defined in  Lemma~\ref{technical}. It is easy to see that the angles $u,w,x,y,z$ satisfy the
assumptions of  Lemma~\ref{technical} (use the fact that $x_2$ is inside the circumcircle
of the triangle $x_1x_3x_4$), so the result is now an immediate consequence of
Lemma~\ref{technical}.
\end{proof}

As a rather simple corollary of the last theorem we get the following result.

\begin{theorem}\label{convex}
Let $x_1, x_2, x_3, x_4$ be vertices of a convex quadrilateral. Then Conjecture~\ref{conj2} holds
for $x_1, x_2, x_3, x_4$.
\end{theorem}

\begin{proof}
We need to prove that the inequality (\ref{formula}) holds. It suffices to show that
\begin{equation}\label{two}
3+\cos\alpha_{l,i}+\cos\alpha_{l,j}+\cos\alpha_{l,k}\geq 2
\end{equation}
for $l=1,2,3,4$. Indeed, then the left hand side of (\ref{formula}) is greater than or equal to
twice the left hand side of (\ref{inequality}), so (\ref{formula}) follows from (\ref{inequality}).

The left hand side of each of the inequalities (\ref{two}) is of the form
$3+\cos \alpha +\cos\beta +\cos(\alpha+\beta)$ with nonnegative  $\alpha$, $\beta$ such that
$\alpha+\beta\leq \pi$. The result follows now from the identity $1+\cos \alpha +\cos\beta +
\cos(\alpha+\beta)=4\cos(\alpha/2)\cos(\beta/2)\cos[(\alpha+\beta)/2]$.
\end{proof}

\begin{remark}
{\rm
\

\begin{enumerate}
\item Inequality (\ref{inequality}) remains true when one of the points $x_1,x_2,x_3,x_4$ is inside
the triangle formed by the remaining three points. This follows from an appropriate version of
Lemma~\ref{technical}, which can be proved along the same lines (basically it is the same lemma
but for $w,x$ which are both negative and with some of the assumptions slightly adjusted).
However, one of the inequalities (\ref{two}) is no longer true in this case so our derivation
of Conjecture~\ref{conj2} is no longer valid.  Nevertheless, the inequality (\ref{inequality}) seems
of independent interest. We have in fact the following conjecture.

\begin{conjecture}\label{conj4}
Inequality (\ref{inequality}) holds for any four distinct points $x_1,x_2,x_2,x_4$ in $\mathbb R^3$.
\end{conjecture}
\noindent
Using R Statistical Software, we have verified this inequality for several million random tetrahedra
so we are quite confident in its validity.

\item Consider any four distinct points $x_1,x_2,x_2,x_4$ in $\mathbb R^3$. Even though the inequalities
(\ref{two}) do not hold in general, it seems that the left hand side of (\ref{formula}) is always greater
than or
equal to twice the left hand side of (\ref{inequality}). Again, we verified this inequality for several
million random tetrahedra so we state it as a conjecture.

\begin{conjecture}\label{conj5}
The left hand side of (\ref{formula}) is greater than or
equal to twice the left hand side of (\ref{inequality}) for any four distinct points $x_1,x_2,x_2,x_4$ in
$\mathbb R^3$.
\end{conjecture}
\noindent
Clearly, Conjectures~\ref{conj4} and \ref{conj5} together imply Conjecture~\ref{conj2}.

\item Formula (\ref{poly}) for a (non-degenerate) tetrahedron contains the term $288V^2$, which
one could hope to incorporate in proving  Conjecture~\ref{conj2}. However, using the result of
Crelle (see Remark~\ref{remcrelle} (3))
one can prove that $288V^2\leq 4d_3(r_{1,2}r_{3,4}, r_{1,3}r_{2,4}, r_{1,4}r_{2,3})$ for any tetrahedron
and the equality holds if and only if the tetrahedron is isosceles (i.e. all its faces are congruent to each
other). In particular, Conjecture~\ref{conj2} is true for vertices of any isosceles tetrahedron.
\end{enumerate}
}
\end{remark}

For the remaining part of this section we will assume that the points $x_1,x_2,x_3,x_4$
are vertices of an inscribed quadrilateral. It follows that the associated Crelle triangle is degenerate
so $d_3(r_{1,2}r_{3,4}, r_{1,3}r_{2,4}, r_{1,4}r_{2,3})=0$ (this is the celebrated Ptolemy's
theorem). Thus Conjecture~\ref{conj2} in this case immediately follows from (\ref{poly}). Our goal is to prove
Conjecture~\ref{conj3} in this case. As noted in \cite{east}, Conjecture~\ref{conj3} can be expressed as
follows:
\begin{equation}\label{con3}
\left|D(x_1,x_2,x_3,x_4)\prod_{1\leq i<j\leq 4}(2r_{i,j})\right|^2\geq \prod_{l=1}^{4}\left( d_3(r_{i,j},
r_{j,k},r_{i,k})+8r_{i,j}r_{j,k}r_{i,k}
\right)
\end{equation}
(recall our convention that $\{i,j,k,l\}=\{1,2,3,4\}$). This is based on a rather simple observation that
$8\cdot AB\cdot AC\cdot BC \cdot D(A,B,C)=d_3(AB,AC,BC)+8\cdot AB\cdot AC\cdot BC$
for any three points $A,B,C$ in $\mathbb R^3$ (see also Remark~\ref{3points}).
When the points $x_1,x_2,x_3,x_4$ are coplanar,
the formulas (\ref{poly}), (\ref{poly1}), (\ref{assoc}) allow us to state (\ref{con3}) in an equivalent
form as follows:
\begin{eqnarray}\label{con31}
\left(16+\sum_{l=1}^{4} (3+\cos\alpha_{l,i}+\cos\alpha_{l,j}+
\cos\alpha_{l,k})\delta(x_ix_jx_k)-2\delta(ABC)\right)^2\geq \\ \nonumber
\prod_{l=1}^{4}\left( \delta(x_ix_jx_k)+4\right),
\end{eqnarray}
where $ABC$ is the associated Crelle triangle.

Our last result is the following theorem.

\begin{theorem}\label{inscribed}
Conjecture~\ref{conj3} holds for the vertices of an inscribed quadrilateral.
\end{theorem}
\begin{proof}
Let $x_1x_2x_3x_4$ be an inscribed (convex) quadrilateral. Since the associated Crelle triangle is degenerate,
by (\ref{con31}) we need to prove that
\[\left(16+\sum_{l=1}^{4} (3+\cos\alpha_{l,i}+\cos\alpha_{l,j}+
\cos\alpha_{l,k})\delta(x_ix_jx_k)\right)^2\geq
\prod_{l=1}^{4}\left( \delta(x_ix_jx_k)+4\right).
\]
Let $A_l=1+\cos\alpha_{l,i}+\cos\alpha_{l,j}+
\cos\alpha_{l,k}$, $B_l=\delta(x_ix_jx_k)=\cos\alpha_{i,l}+\cos\alpha_{j,l}+
\cos\alpha_{k,l}-1$ for $l=1,2,3,4$. Thus we have to prove that
\begin{equation}
\left(16+\sum_{l=1}^{4} (2+A_l)B_l\right)^2\geq \prod_{l=1}^{4}\left(B_l+4\right).
\end{equation}
Since the quadrilateral is inscribed, it is easy to see that $A_1+A_3=A_2+A_4=B_1+B_3+4=B_2+B_4+4$
and $A_l-B_l=2+2\cos(\angle x_{l-1}x_lx_{l+1})$ for $l=1,2,3,4$. In particular, $A_l\geq B_l\geq0$ for
$l=1,2,3,4$. Note that
\[ \sqrt{(B_1+4)(B_3+4)}\leq 4+\frac{B_1+B_3}{2} \ \text{and}\ \sqrt{(B_2+4)(B_4+4)}\leq 4+\frac{B_2+B_4}{2}. \]
It suffices then to prove that
\[ 16+\sum_{l=1}^{4} (2+B_l)B_l\geq \left(4+\frac{B_1+B_3}{2} \right)\left(4+\frac{B_2+B_4}{2} \right).\]
Since $B_1+B_3=B_2+B_4$, it is enough to show that
\begin{equation}\label{fin}
 8+ 2(x+y) +x^2+y^2\geq \frac{1}{2}\left(4+\frac{x+y}{2} \right)^2
 \end{equation}
holds for $x=B_1$, $y=B_3$ and for $x=B_2$, $y=B_4$. As a matter of fact, (\ref{fin}) holds
for any real numbers $x,y$ as it is easily seen to be equivalent to
\[6x^2+6y^2+(x-y)^2\geq 0.\]
\end{proof}

\begin{remark}
{\rm Inequality (\ref{con31}) is equivalent to Conjecture~\ref{conj3} only for four coplanar points.
In general, for arbitrary four points in $\mathbb R^3$ it only implies Conjecture~\ref{conj3} (i.e. it is a
stronger inequality). Nevertheless, numerical investigation leads us to believe that the following
should be true.

\begin{conjecture}\label{conj6}
Inequality (\ref{con31}) holds for any four distinct points in $\mathbb R^3$.
\end{conjecture}

}
\end{remark}

\vspace{5mm} \noindent {\bf Acknowledgements.} We are grateful to the referee for making thorough
remarks improving our paper.

\end{document}